# Encoding discrete quantum algebras in a hierarchy of binary words

**T. E. Raptis**[abc]

[a]National Center for Science and Research "Demokritos", Division of Applied Technologies, Computational Applications Group, Athens, Greece.
[b]University of Athens, Department of Chemistry, Laboratory of Physical Chemistry, Athens, Greece
[c]University of Peloponnese, Informatics and Telecommunications Dept., Tripolis, Greece
email: rtheo@dat.demokritos.gr

**Abstract.** It is shown how to endow a hierarchy of sets of binary patterns with the structure of an abstract, normed C*-algebra. In the course we also recover an intermediate connection with the words of a Dyck language and Tempereley-Lieb algebras for which we also find that an effective arithmetic code is possible albeit of greater complexity. We also discuss possible applications associated with signal theory and waveform engineering on possible ways to embed discrete computational structures in an analog continuum substrate.

## 1. Introduction

The notion of C*-algebras, with "*C*" originally standing for "closed", originates in the work of Gelfand and Neimark [1], and later Segal [2] as an attempt to an abstract generalization of the work of Jordan and Von Neumann [3] on quantum algebras of operators associated with Heisenberg matrix mechanics. Recent advances in radar waveform synthesis [4] and [5] reintroduced the theme of Heisenberg and other group structures as important tools in the general signal processing theory. As a matter of fact, an important observation is already met at page 148, in the last chapter of [5] directly associating signal theory with an older paper by Birkhoff and Von Neumann [6] on the foundations of quantum logic.

Recent attempts in understanding the problem posed by the famous Bekenstein entropic constraint [7], [8] implying a finite information limit led to various attempts for remedying the situation [9], [10], [11]. Together with the notion of finitary information structures as being fundamental in the understanding of quantum dynamics seem to imply a tighter connection between the two themes of quantum structures and general signal theory. In [12], the author attempted the introduction of a method for embedding certain automata and computational structures in a system of continuous signals. Additionally, one can notice the existence of some dispersed evidence in the literature for the existence of certain analog models of computation that appear to be of similar power with models of quantum computation (indicative references in [13] – [18]).

Notably, in certain distributed or holographic models of computation, one also finds elements of non-locality as in the work of Smolenski [19] and Smolenski and Fernandez [20] where Kronecker products are used as abstract knowledge representations or in the work of Plate [21], [22] where superpositions of convolution products are used for similar purposes. As a result of this original research, Aerts [23] and later Aerts and Gabora [24] showed the equivalence of such models with the State Context Property (SCOP) formalism, a generalization of quantum mechanics formalism [25], [26].

Complementary to these evolutions, Aravind [27] and later, Kauffman [28] had already proposed a direct analogy between quantum and topological entanglement (linking) which led to finding many intriguing connection between knots [29], [30] the Artin's Braid group [31], [32] and the associated Temperley-Lieb (TL) algebras [33], [34]. At about the same time, a theorem due to Brylinski [35] allowed the complete solution of the transcribed relations of the set of TL generators into tensorial product representations known as the algebraic Yang-Baxter equation [36], [37]. Based on this, a type of universal gates for topological quantum computing was found [38], [39], [40], [41], [42] of which the analog on a signal theory is of distinct interest and shall be explored in the following sections.

Attempting to unify these diverse themes under the umbrella of general signal theory and signal processing imposes a question on how much of the underlying methods and axioms can actually be "emulated". The particular term should be used with caution since, any two theoretical constructs capable of mimicking certain effects of each other could be said (for instance, "analog gravity"[43]) to be emulators of one another and yet, if at a certain point the emulation becomes so faithful that no further questioning of each system can reveal any difference then the two have to be considered as identical or isomorphic, even if the two appear to be based on a different ontology.

In the next sections, we attempt to first set up a unique framework allowing the examination of several different themes under the common theme of their arithmetized encoding, a kind of Gödelianization method for bounded string tuples. In section *2*, we prove that a direct transfer of the axiomatic framework of abstract $C^*$-algebras is possible for a hierarchy of lexicographically ordered (LEX) dictionaries of binary words using the construct of an Inductive Combinatorial Hierarchy (ICH) that was introduced in [44] as a generic toolbox for the study of global maps of automata. In section *3*, we observe that the result of the necessary conjugation introduced in the previous section, sends directly to a special subset of words of a Dyck language [45], [46] which is isomorphic with the diagrammatic interpretation of the elements of any $TL_n$ algebra. We then examine the closure of any $TL_n$ product inside the ICH and we attempt the complete arithmetization of the resulting subset. In section *4*, we examine the possible translation of the previous constructs in the Fourier domain of a system of signals that could be processed by a class of, spectral analog machines based on special convolution filtering. We close with a discussion of certain ambiguities involved in signal observation that we believe, they are also closely tied with the continuous versus discrete dichotomy, as well as the need for finite information structures on a continuous background.

## 2. A hierarchy of $C^*$-algebras

An Inductive Combinatorial Hierarchy is defined on an inclusive sequence of integer intervals $S_0 \subset S_1 \subset ... \subset S_k \subset ...$ with $S_0 = \{0,1\}$, $S_k = \{0,...,2^k - 1\}$ associated with a self-similar hierarchy of $k$ x $2^k$ Boolean matrices or "word dictionaries" in LEX order, $W_0 \subset ... \subset W_k \subset ...$ characterized by the fact that each serves as a complete dictionary of all possible patterns when each interval ends with a Mersenne number. Let then an arbitrary automaton or algorithm performing a mapping $M : S_n \to S_m$ which can always be made an endomophism in $S_{Max(n,m)}$. We first define a pair of pull-back and push-forward maps as a "*Decoder*" $d : N \to W$ and an "*Encoder*" $e : W \to N$, satisfying $ed = de = I$ so as to form a chain $S_k \xrightarrow{d} W_n \xrightarrow{M} W_m \xrightarrow{e} S_k$. The global map of $M$ allows forming arithmetized codes as $f : N \to N : f = eMd$. Moreover, we can interpret $M$ as a "parsing" map in higher radices $b > 2$ for the associated dictionaries $W_i^b$, such that given a set of generators of some algebra $e_i$, $i = 1,\ldots,b$ and a functional lookup table (FLUT) of exactly $b$ indexed functions or algorithms, $M$ sequentially assigns all indices from any word in $W_i^b$ to a set of members of the particular FLUT deriving the particular sequence of actions of any word $e_ie_j$... The result over any $W_i^b$ may be called, the hierarchy of global maps for each such structure and it is an empirically conjectured result that more often than not, the resulting global sets inherit a fractal structure which originates in the inner combinatorial structure of $W_i^b$ per se, apart from the details of any particular algebraic symmetry or group structure. This hierarchy then forms the complete product of all possible combinatorial structures or combinatorial designs [47],[48] given additional filters as indicators over any $S_k$, with all possible algebras.

We then ask to construct the equivalent of some discrete normed division algebra over binary words of constant length in any $W_k$ via actions on their integer values associated via the polynomial representation over a complex field ***C***. This calls for the separate definition of an addition and multiplication for word patterns plus the existence of an involutive conjugation operation, usually denoted as (*) and an appropriate norm $\|...\|$ satisfying the C*-identity, to guarantee the uniqueness of the norm as $\|xx^*\| = \|x\|^2$. To distinguish between an integer representation of a word and an ordinary number in ***C***, we shall denote the former as $(\nu)_2 = \sigma_1\sigma_2...\sigma_k, \sigma_i \in \{0,1\}$. We shall further interpret the conjugation operation in its standard meaning for complex numbers while for binary representations we shall adopt the equivalent of the action of two involutive operators, one given as a bitwise complement operator or arithmetically, $\bar{\nu} = 2^k - \nu - 1, \ \forall \nu \in S_k$ and a reflection operation which inverts the order of bits in any pattern associated with $v$ as a left-right flip and extracting a new integer which we shall denote as $v' = \sim v$. The two operations naturally commute and their combined operation shall be denoted as the conjugate pattern represented as $\nu^* = \sim \bar{\nu}$. As a composite of two involutions, conjugation is a bijective automorphism in any $S_k$. The particular operation is a "folding" with respect to the coding map, as given by the polynomial representation which is depicted schematically in figure 1, with the arrows denoting the direction of the $2^n$ coefficients.

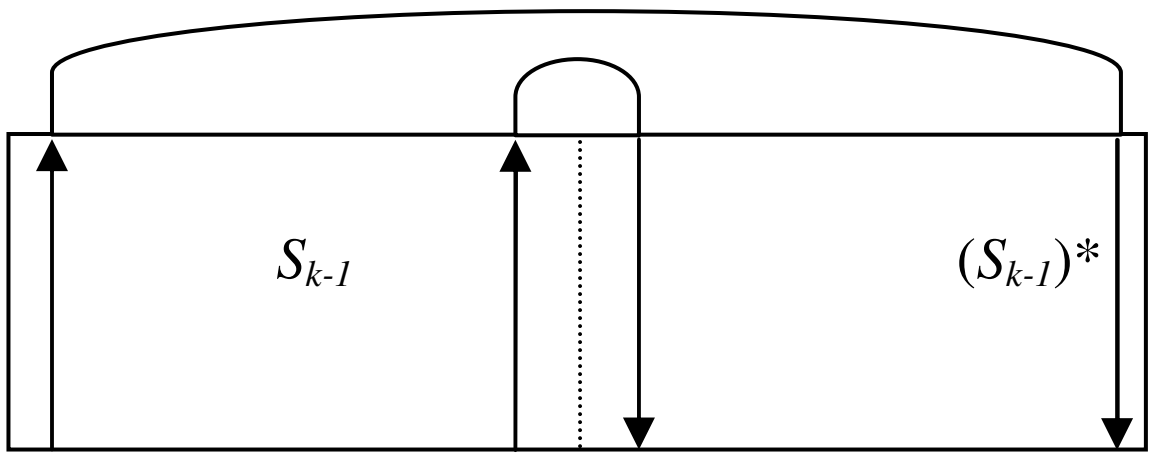

**Fig.1** Schematic depiction of the action of the weighted average over all conjugate pairs over the ordered array of all words in a total interval $S_k$. with arrows showing the $vv^*$ product while reverse arrows hold for $v^*v$.

From the independent action of the two involutions of complementation and reflection we also obtain the identity

$$(\sim \nu)_2 + \sim \left(2^k - (\nu)_2 - 1\right) = \nu + \bar{\nu} = 2^k - 1 \qquad (1)$$

We notice that reflections admit a fixed point subset of palindromes in any $S_k$ for which $\sim \nu = \nu$. Evidently, for any even order $S_{2\lambda}$ there will be exactly $2^\lambda$ palindromes while for odd orders $2\lambda + 1$ one must stabilize the middle bit thus leaving again $2^\lambda$ fixed points. The previous definitions allow the standard addition as

$$\left(z_1\mu + z_2\nu\right)^* = z_1{}^* \mu^* + z_2{}^* \nu^*, \ \ \forall z_1, z_2 \in \mathbf{C}, \ \ \forall \mu, \nu \in N \qquad (2)$$

For the case of multiplication we shall introduce the non-commuting string concatenation as

$$\mu\nu = \left[(\mu)_2,(\nu)_2\right] = (\mu + 2^k \nu)_2 : S_k \to S_{2k}, \forall(\mu)_2,(\nu)_2 \in S_k \tag{3}$$

We notice that (3) defines an endomorphism from both even-odd orders to an even order word set. We must also supplement each set with the null or empty string "{}" as a neutral element. Then we can define an inverse via the complements such that $\nu\bar{\nu} = \nu\nu^{-1} = \{\}$. We notice that complementation entails a natural span of any interval into inverse subsets as $S_k = S_{k-1} \cup S_{k-1}^{-1}$. It holds by construction that conjugation is distributive such that

$$(\mu\nu)^* = \left[(\nu)_2{}^*,(\mu)_2{}^*\right] = \nu^* + 2^k \mu^* = \nu^* \mu^* \tag{4}$$

We can also find invariants of the above through the symmetric form $z(\mu\nu^*) + z^*(\nu\mu^*)$ as well as the symmetric product $\nu\mu^*\mu\nu^*$.

From the multiplication definition it becomes apparent that in order to satisfy the norm axiom for words we shall need to take exponents of some additive map over the integers in any $S_k$. Assume then a word norm as

$$\left\|(\nu)_2\right\| = e^{gf(\nu)} \tag{6}$$

with the choice $g = ln(2)$ being a special case. We immediately arrive at the simultaneous conditions

$$\begin{aligned} \|\nu^*\nu\| &= \left\|(\nu^* + 2^k\nu)_2\right\| = \|\nu\nu^*\| = 2^{f(\nu^*+\nu)} \\ \|\nu\|\|\nu^*\| &= 2^{f(\nu^*)+f(\nu)} = 2^{2f(\nu)} \end{aligned} \tag{7}$$

The above leads to a variant of the Cauchy functional equation [49], [50] of which there are many non-trivial solutions. We are particularly interested in finding some method by which $f\left((\nu)_2^*\right) + f\left((\nu)_2\right) = f\left((\nu^* + \nu)_2\right)$ plus the symmetry condition $f\left((\nu)_2^*\right) = f\left((\nu)_2\right)$ where the details of each particular pattern in $(\nu)_2$ should be taken into account by the particular choice of $f$.

An obvious solution for binary words would be given with the aid of an important function in many areas like coding theory and number theory [51] namely the digit-sum $\Sigma_2(\nu)$ which satisfies

$$\Sigma_2(\nu + \mu) = \Sigma_2(\nu) + \Sigma_2(\mu),\ \ \Sigma_2(\nu^*) + \Sigma_2(\nu) = \Sigma_2(\nu^* + \nu) = k \tag{8}$$

While this choice would trivially satisfy the last $C^*$-axiom, it is still highly degenerate since all string combinations in $S_k$ x $S_k$ would share the same norm value $2^k$. We notice that by including the empty string, this could lead to a particular finite, unital hierarchy of algebras by a redefinition of the norm as $2^{-k} e^{gf(\nu)}$. There is certain evidence in the literature that if stable, such can be embedded in a so called AF-algebra [52], [53], [54] as an inductive limit but the topic is beyond our specific interests here.

To remedy the degenerate structure of the previous choice, we need to find another method of characterization of all patterns by a sufficiently sensitive measure of each pattern's complexity and structure. As is evident, such cannot be given in terms of the binary Shannon's entropy since for any finite string, the particular function depends solely on a probability defined as $\Sigma_2(\nu)/k$ thus giving the same output for each binomial set of $\binom{k}{\Sigma_2(\nu)}$ elements. The only possible way to fix this seems to be the promotion of *f*(*v*) to a full square form with the kernel matrix being that of an appropriate discrete transform over each bit pattern. We also recall that the reflective symmetry of the conjugation operation has similar properties with the symmetry $f(-\omega) = f(\omega)$ in a Fourier domain of a symmetric function. Indeed, taking the summands of the real cosine (DCT) and sine (DST) transform can solve this problem as shown in figure 2.

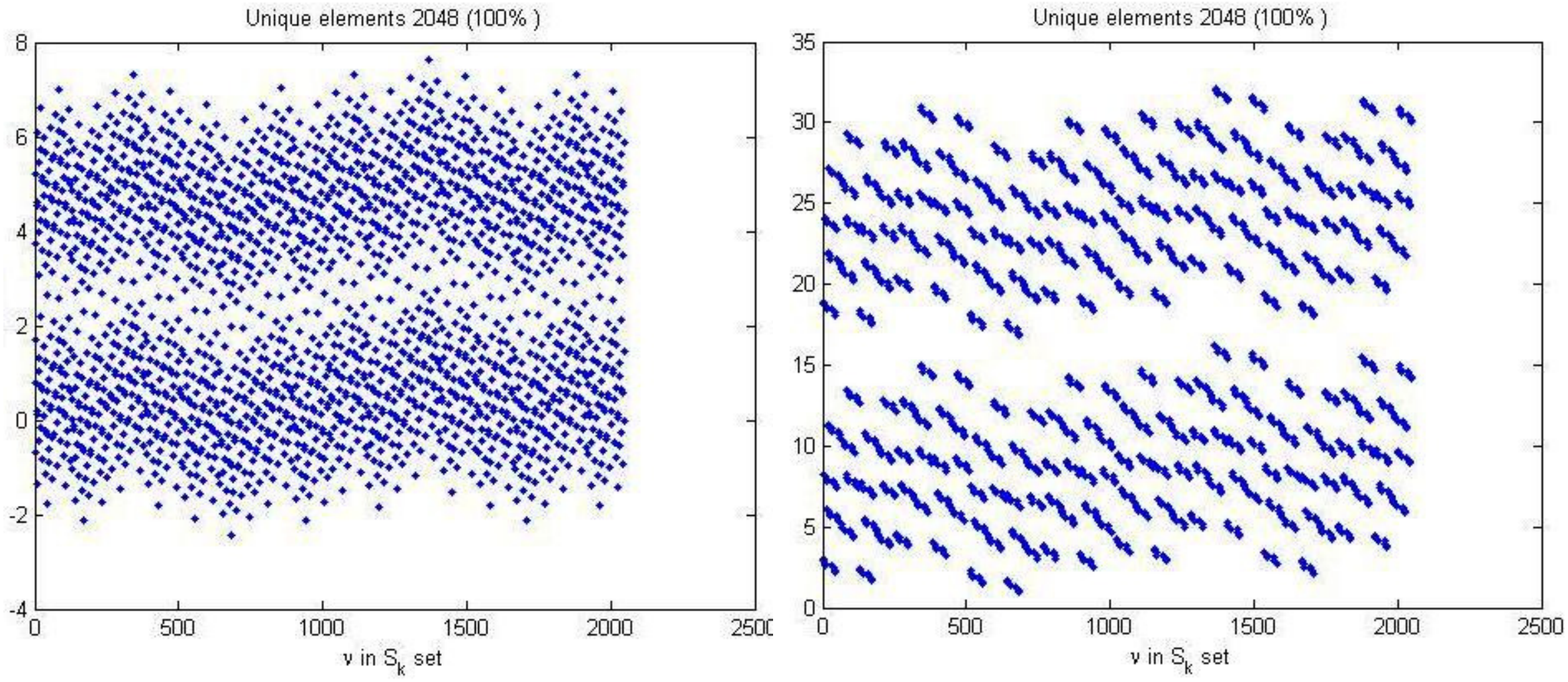

**Fig. 2** The resulting log-norms for all 11-bit words in $S_{11}$ with (a) the DCT and (b) the DST showing full degeneracy lifting.

We can then redefine (6) with $F_{C,S}$ denoting the DCT or DST respectively,as

$$\left\|(\nu)_2\right\| = e^{g\sum_{i=1}^{k} F_{C,S}(\nu)_2} \tag{8}$$

The above representation is not unique and other types of transforms like the discrete auto-correlation function, also associated with the DFT of each pattern could be used. In the next section we examine some interesting relations between the conjugation operation and the alternative theme of the Tempereley-Lieb algebras.

### 3. Dyck words and Temperley-Lieb algebras

Let $D_k \subset S_k$ the subset of words of a Dyck language of length *k* with equidistributed zeros and ones. The number of words across the hierarchy is known to increase as the sequence of Catalan numbers $(k+1)^{-1}\binom{2k}{k} = (k+2)...(2k)/k!$. It is possible to extract the particular Dyck subsets over a hierarchy of total word sets using an arithmetic encoding of the form "(" → "1" and ")" → "0" and vice versa. This result in two complementary isomorphic sets since any $S_k$ is characterized by a 2-span of dual subspaces and any Dyck word must either start with 1 and end with 0 or vice versa for the complementary interval $\overline{S}_{k-1}$.

In order to properly apply the total set of constraints upon all the words of any $S_k$ we need an intermediate tool for the characterization of the interior bit block structure of any word. This can be given via a special map, establishing a lossless compression scheme in an alphabet of higher radix $k$ as follows. Let a transcription map from the original binary words to the words of a Hadamard code {-1, +1} via $\sigma_i \to s_i = 2\sigma_i - 1$ where $(v)_2 = \sigma_0\sigma_1...\sigma_{K-1}$. Let also $\boldsymbol{c}$ denote a vector of coefficients of a block counting polynomial and a map from each word to a distinct associated polynomial defined as

$$(v)_2 \to P_2^C(v) = \sum_{i=1}^{\delta(v)} c_i z^i : c_i = \begin{cases} sign(s_i), & d = 0 \\ dc_i + sign(s_i), & d \neq 0 \end{cases} \quad d = \frac{1}{2}\left(s_i sign(s_{i-1}) + 1\right) \tag{9}$$

The above set of polynomials has several nice properties that will be studied elsewhere. The resulting set of coefficients for any polynomial is an alternating sequence with the sign following the sequence (*mod 2*) and it satisfies by construction

$$\sum_{i=1}^{\delta(v)} | c_i | = k \tag{10}$$

Due to (10) and the fact that any $S_k$ is the complete combinatorial superset of order $k$ for every $\{0,1\}^k$ unit hypercube, the coefficients set necessarily realize all integer partitions of $k$ although repetitively and not in a lexicographic order. The varying degree $\delta(v)$ shall be called the block dimension and its evolution across the hierarchy reveals a fractal sequence related to the paperfolding or dragon-curve sequence. We shall only give without proof an exact, recursive generator across any maximal interval $S_k$ as

$$\delta_{k\max} \leftarrow \left[\delta_{k-1}, \sim(\delta_{k-1}) + d\right], \;\; d = 1 \;\; iff \;\; k < k\max \tag{11}$$

Recursive relations like (11) can often be found using successive sequence folding (reshaping into matrix form) until the simplest possible reproducing maps can be located as explained in [44]. Notably, the class of polynomials defined in (9) can be extended into arbitrary higher radices $b > 2$ using complex coefficients as $| c_i^j | \omega_0^j, j = 0,...,b-1$ where $\omega_0$ the primitive $b$th root of unity, but they are of no concern here.

The Dyck constraints on the partial summands can now be translated as

$$\sum_{i=1}^{j} c_i > 0, \;\; j = 1,...,\delta(v) - 1, \;\; \sum_{i=1}^{\delta(v)} c_i = 0, \;\; \delta(v) = 2\kappa \tag{12}$$

The last condition on the block dimension being even is a natural consequence of the equidistribution property which must also hold for symbol blocks. Regarding the conjugation operation (*) we observe that no block merging can happen since all words formed this way result in symbol sequences having the form $\sigma_0,...,\sigma_{k-1}, 1-\sigma_{k-1},...,1-\sigma_0$. The resulting block coefficients will necessarily satisfy $\boldsymbol{c}^*$ = -(~$\boldsymbol{c}$) such that $\mathbf{c} + \mathbf{c}^* = \mathbf{0}$. If a particular $\boldsymbol{c}$ vector satisfies the Dyck conditions (12) so will do the total [$\boldsymbol{c}$, $\boldsymbol{c}^*$] since the concatenation of two Dyck words is a Dyck word and the simultaneous sign change followed by order reversal in the summands does not affect the positivity of the partial summands. The total Dyck set then admits a natural split as $D_k = D_{k-1} \cup D_{k-1}^*$.

From the definition of multiplication (3) and from (7) all concatenation products $S_k \times S_k^*$ related to the norms of section 2, we can define an endomorphism into a subset of $D_{k+1}$ as $D_k \xrightarrow{v+2^k v^*} D_{k+1}$. It is also known that any $D_{k+1}$ coincide exactly with the number of words from a Temperley-Lieb or $TL_n$ algebra of which the generators satisfy the Jones relations. There are also certain homomorphisms associating the generators of a special type of Hecke algebra with those of a $TL_n$ algebra as well as those of the Artin's Braid group. We conclude that there is always a cross section $C^* \cap TL_n$ between a full dictionary of words endowed with the $C^*$-algebra of section 2, and a subset of TLn words.

We are particularly interested in an isomorphic, diagrammatic representation of a $TL_n$ which can be given straightforwardly with an input-output model of an enclosing box. Each such box has the same number of inputs and outputs of which the connectivity is given via a pair of objects we shall term the "connectors", being either straight or diagonal line segments between the upper inputs and lower outputs and the "cups" or hemi-circles, acting like short circuits between either inputs or outputs. The isomorphism with Dyck words comes from the folding of an even bit pattern such that any pair of closed parentheses 1…0 beyond $k$ bits has to be folded thus giving rise to a connector, otherwise it becomes a cup. Additionally, it allows nested cups but disallows crossings between inputs or outputs of cups and connectors. Such a situation is depicted in figure 3 where we also show the exact correspondence with the associated parenthesized Dyck word and its arithmetic codes.

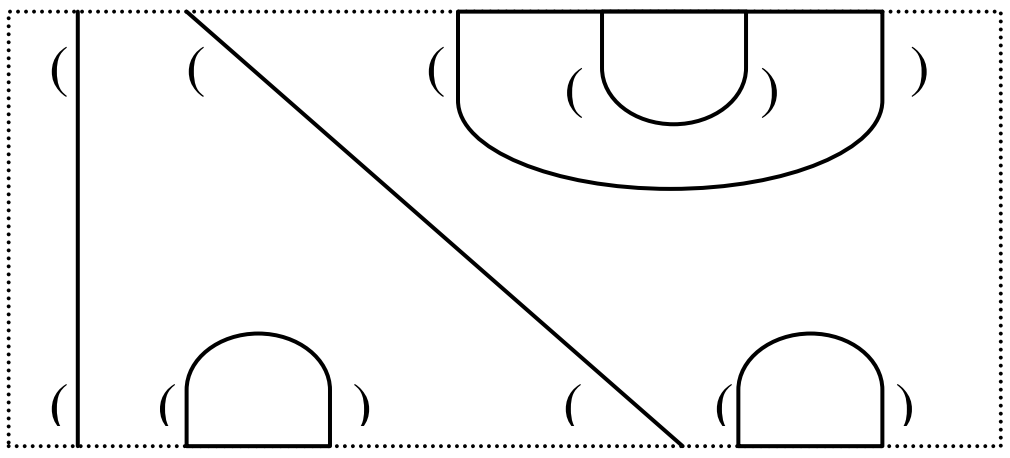

**Fig. 3** Depiction of a diagrammatic $TL_n$ box with its associated Dyck word corresponding to the pair of integer codes, 111100100100 and its complement. The enclosing box is indicated with interrupted lines.

We can then present any $TL_n$ in an arithmetized form. We notice that the particular case of figure 3 falls outside of the subset of conjugate words of the previous section since the total word does not always have the symmetry of a ($v$, $v^*$) pair. Products in any such diagrammatic construct correspond to a superposition of boxes in which a special simplification algorithm is to be applied. Specifically, all cups in the interface between two boxes act as annihilators either for input elements of the lower box, or for output elements of the upper box. They annihilate either a pair of connectors or another cup when they fully overlap with an input-output relation. They also partly annihilate any input output chain of partially overlapping connectors and cups leaving either single cups or a diagonal connector. The two boxes then collapse into one. The total operation is a mapping $D_{k+1} \times D_{k+1} \to D_{k+1}$ so that the set of all Dyck words in $S_{k+1}$ becomes a closure for this product at any level of the hierarchy. In figure 4(a), we show the subset of all Dyck words in their arithmetic representations inside the total $S_{10}$ interval of all 10-bit words. A mirror image exists for the complementary encoding in the other half-interval. In 4(b) we also performed the analysis of the associated Boolean indicator over the first half in $S_9$, by finding the coefficients of its block counting polynomial to make evident the inherent self-similarity.

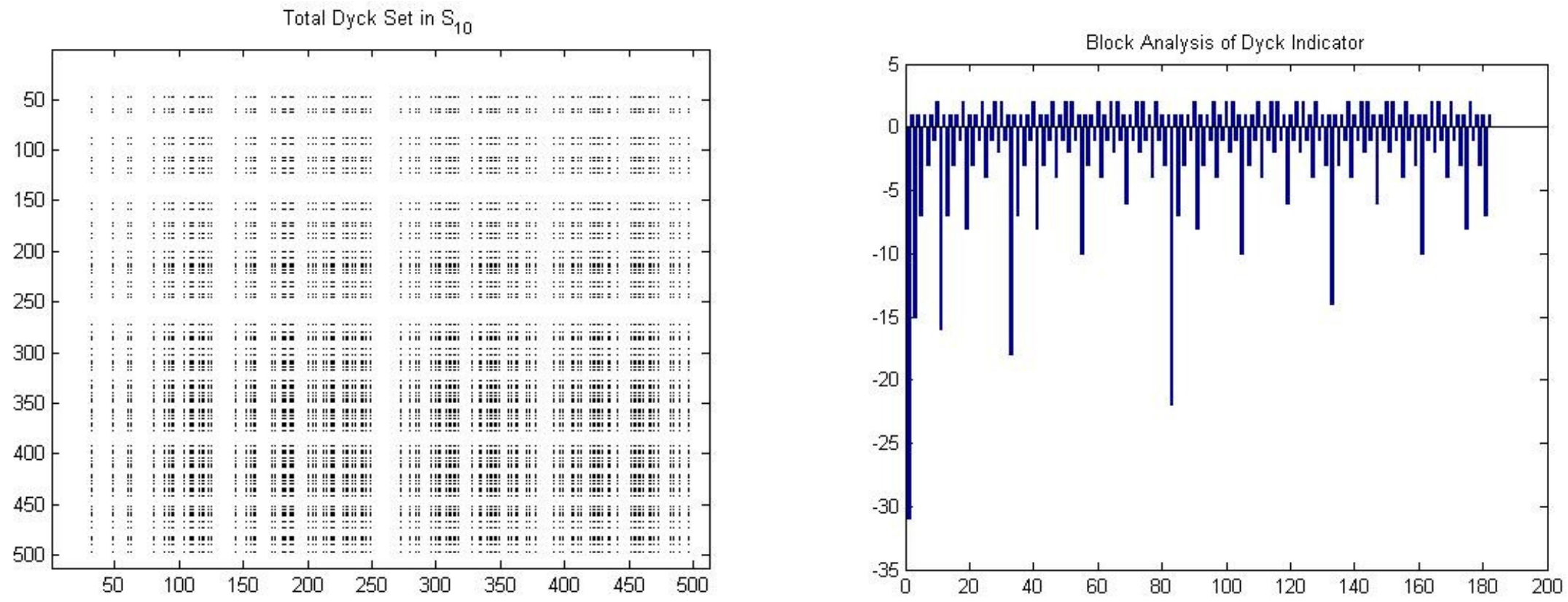

**Fig. 4** (a) The $2^9$ x $2^9$ matrix of all valid $TL_8$ arithmetized product factors for the associated Dyck words set in $S_{10}$ with reverse colormap ( black = 1, white = 0) and (b) the coefficients of the block polynomial for the associated Boolean indicator over all $S_9$.

In what follows we shall try to extract an algorithmic simplification of lowest possible complexity for all these products using their arithmetic encoding which would be equivalent into a logical circuit or a composition of bitwise Boolean functions for the arithmetized factors. To this aim, we define an abstract product between the associated Dyck words as

$$(\nu \circ \mu)_{TL}^{k} = f(\sim \nu, \mu) = f(\nu, \sim \mu), \quad \forall \nu, \mu \in D_{k+1} \tag{13}$$

In (13), the two items on the *rhs* show the equivalence in matching expansion symbols in reverse order for each of the terms since, any integer expansions have to be overlapped with the second half of $(\nu)_2$ over the first half of $(\mu)_2$ and vice versa, in order to fit the diagrammatic "box" picture. The question then is how to find an equivalent algorithm or rules of an automaton isomorphic to the final action of the *f* map.

We observe that there are two major classes of products of which the first is characterized by the fact that the final product is composed directly from the half brackets of the upper and lower boxes independently of the presence of annihilation acts in their middle interface. In any other case, information is transferred from the middle interface (bit block alterations) to one or both of the final input-output pairs. An example is given in figure 5.

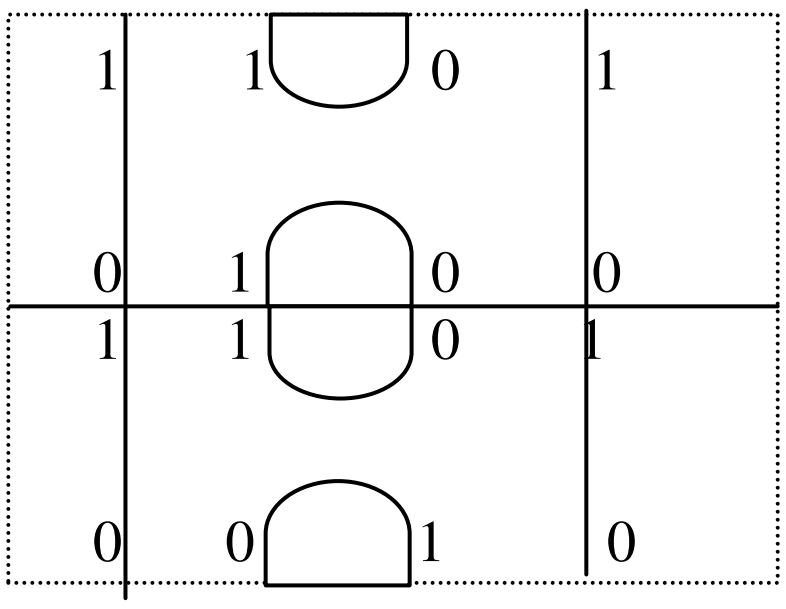

**Fig. 5** An example of a TL product in which the upper input and lower output structure of the enclosed fluxes remain invariant.

We can then immediately reduce (14) assuming a Boolean function *g* deciding between the two classes, using the following instrumentation. Any pair shall be written as $\nu = \nu_u + 2^k \nu_d$, $\mu = \mu_u + 2^k \mu_d$ and we shall use $h_u(n) = \mathrm{mod}(n, 2^k)$ and $h_d(n) = \lfloor n2^{-k} \rfloor$ as sub-word extractors so as to be able to write

$$\nu \circ \mu = \begin{cases} \nu_u + 2^k \mu_d, & g(\nu_d, \mu_u) = 1 \\ f(\nu_u, \sim \nu_d, \mu_u, \mu_d), & g(\nu_d, \mu_u) = 0 \end{cases}, \quad \forall \nu, \mu \in D_{k+1} \tag{14}$$

For all other cases, a special simplification algorithm is required to alter the block structure associated with bracket rearrangements that are mostly non-local. In (14), the $f$ map can also be written as a transition rule of the form

$$f(\nu_u, \sim \nu_d, \mu_u, \mu_d) = \begin{cases} \nu_u \xrightarrow{R(\sim \nu_d, \mu_u)} \nu'_u \\ \mu_d \xrightarrow{R(\sim \nu_d, \mu_u)} \mu'_d \end{cases} \tag{15}$$

The expression $R(\sim \nu_d, \mu_u)$ is used to denote the rest of the information consumed at the interface. Due to the restrictions of any Dyck language, any information transfer from the interface can be analyzed into a set of transitions of previously open brackets to closed ones inverting a single bit in each such case.

An important difference between the two branches of (15) is that each consists of complementary transitions. Specifically, the first branch can only alter two open brackets into a closed pair hence a transition 1 → 0 while in the second branch, significance of brackets gets reversed hence only transitions 0 → 1 can happen. This allow further simplifying the second branch of (15) assuming *R* returns a logical mask $r_u + 2^k r_d$ such that one can simultaneously achieve the reversal effect in all relevant positions via

$$f(\nu_u, \sim \nu_d, \mu_u, \mu_d) = (\nu_u - r_u) + 2^k (\mu_d + r_d) \tag{16}$$

We remind that connector pairs can only be annihilated giving rise to new cups. We notice that the oscillating character of the dynamics in (16) if applied recursively, can lead at particular instances of the general class of batrachion sequences [55], [56].

From the composition laws for the generators $T_i$ of $\mathrm{TL}_k$ we can extract a case of invariance for any $(\nu_u, \mu_d)$ pair. Since for all $i < n$, we have $T_i^2 = cT_i$, this must correspond to a null transition for (16). Moreover, any generator is characterized from $\nu_u = \nu_d$, $\mu_u = \mu_d$. It is elementary to verify that any other composition out of circular permutations of the form $T_i \, T_{\mathrm{mod}(i+j,k)}$ leads to the $(\nu_u, \mu_d)$ pair being invariant for any even number of steps while for odd steps we have an alternating $\nu_d \to \nu'_d = \nu_d - 2^i \quad \mu_d \to \mu'_d = \mu_d + 2^i$ transition where *i* mark a position left of an upper or lower cup respectively.

In any other case as those containing nested cups, there can be longer, non-local annihilating paths which require the *R* map to incorporate a pattern recognition method. In figure 6, we present some limited examples where it is shown that many annihilating patterns are composed of conjugate pairs ($\nu$, $\nu$*) due to the symmetry requirements posed by the closed cup structure.

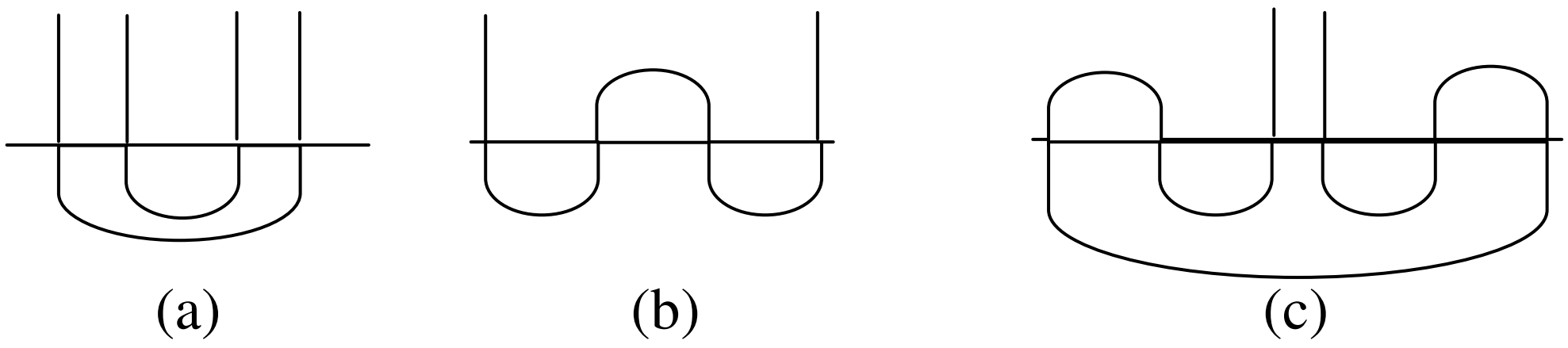

**Fig. 6** Examples of long annihilating paths where the interface values of $\mu_u$ are conjugated pairs. In (a) we have $v_d = 1111$, $\mu_u = 1100$, (b) $v_d = 1101$, $\mu_u = 1010$, (c) $v_d = 101110$, $\mu_u = 110100$.

A method for dealing with such paths can be based on a compressed representation of any superposition of $\begin{pmatrix} v_d \\ \mu_u \end{pmatrix}$ pairs at the interface into a quaternary alphabet as $\left(v_d * +2^k \mu_u\right)_4$. The conjugation of the upper part is used to turn brackets in the reflected sequence so as to ease the comparison with the lower part with brackets on equal footing. Any long, annihilating path then starts always with a superposition as $\begin{pmatrix} 1 \\ 1 \end{pmatrix}$ and ends exactly before the next one. An additional mechanism is required for deciding if one or both of ($v_u$, $\mu_d$) have to be altered as well as the relevant bit pointers for the output Boolean mask in (17). Effort is underway to reduce the complexity of the resulting code into the simplest possible circuit before examining its analog transcription as explained in the next section where we discuss possible applications of the previous sections in signal theory for the purpose of building certain variations of spectral, analog computing machines.

**4. Possible translations to a signal theory**

A spectral based analog computing machine can be defined abstractly as operating on a set of signals $s_i(t)$ of which the frequency or wavelength domain denoted as $\tilde{s}_i(\omega)$ or $\tilde{s}_i(\mathbf{k})$ contains an appropriate encoding of a set of words in a symbolic alphabet. Due to the existence and general properties of Fourier multipliers [57], [58] of which the properties are generalized in the case of arbitrary convolution filters, it is guaranteed that many operations involving word transitions can be realized via an appropriate selection of convolution filters of which the frequency domain application will directly correspond to either $\tilde{s}_1(\omega)\tilde{s}_2(\omega)$ for a pair of signals or $\tilde{K}(\omega)\tilde{s}(\omega)$ for the alteration ("transition") of a particular signal content with $K(t)$ an appropriate kernel function.

Moreover, there do exist certain reduced architectures also known as one-instruction-set computers (OISC) [59], [60], [61] in which a single function also termed "subleq" or "addleq", composed of a single adder with a branching instruction suffices for any other high-level operation via a functional composition with a varying parametrization. An adder is always realizable via a recursive combination of bitwise XOR and AND operations which can also be parallelized as in a symmetric two rule cellular automaton [62] while its transfer into a signal background becomes possible with additional logical masks for the transitions as in (16) thus making any such bitwise operation strictly additive via $\tilde{s}(\omega) + \tilde{s}_{0->1}(\omega) - \tilde{s}_{1->0}(\omega)$ operations. Using similar techniques, it is possible to find distributed holographic representations of elementary systems like cellular automata for which a full multi-linear expression vi a linear filter cascade has been reached including such universal rules as the Turing complete rule 110, which will presented in a subsequent publication [63].

In [12] we presented a particular coding method for waveform synthesis based on the use of what is known as a "Manchester code"[] ("0" → "10", "1" → "01"),a limited form of conjugate pairs ($v$, $v$*), directly into the frequency or the wavelength domain so as to keep also the total spectral power constant apart from localized fluctuations during filtering. Additionally, it is possible to use a kind of hyperspherical code (we use the term to discriminate from existing use of "spherical codes" which belong to a different realm of applications.) Specifically, it is possible to use αν alternative form of the block counting polynomials with positive coefficients using an additional term $c_0$ for encoding the sign resulting in degree $d(v)+1$ polynomials as $mod(v,2) + zP(v)$ in any $S_k$ to define a set of signals as

$$s(t) = \sum_{i=1}^{d(v)} w_i e^{\mathbf{i} n\omega_0 t}, \quad w_i = \sqrt{\frac{|c_i|}{d(v)+1}} : \sum_{i=1}^{d(v)} |w_i|^2 = 1 \tag{18}$$

In (18) we also use a well known technique called orthogonal frequency division multiplexing (OFDM) [64], [65], [66] currently in use for communication applications in which a set of harmonics leads to a kind of Hilbert space representation, with the orthogonality represented by a temporal average over products of $\exp(\mathbf{i} n\omega_0 t)$ terms.

Next, we can ask whether the generator relations for the $TL_k$ algebras could be transcribed in a system of signals $\{s_i^e(t)\}$. The set of composition relations is difficult to interpret in the context of signal convolution since commutativity is not fully allowed since we should have $s_i \otimes s_j = s_j \otimes s_i$ only for $|i-j|>1$. The interest in TLn algebras here stems from their connections with the generators of the braid group $B_n$ satisfying

$$\begin{cases} \sigma_i \sigma_j = \sigma_j \sigma_i, & |i-j|>1 \\ \sigma_i \sigma_{i+1} \sigma_i = \sigma_{i+1} \sigma_i \sigma_{i+1}, & \forall i \end{cases} \tag{19}$$

In the representation theory of $B_n$, the second relation takes the form of the algebraic Yang-Baxter equation of which the solutions were first prescribed in the work of Brylinski [35]. This was followed by intense research on a type of universal gate that could be used for topological quantum computing.

We can examine an alternative possibility for a phase space realization of this representation via the most general Linear Canonical Transform (LCT) which summarizes a lot of other useful transforms that can also be implemented with appropriate optical setups. The most general expression for a unitary LCT is given as

$$\tilde{s}(u) = O_F^{\mathbf{R}}(s(t)) = (-\mathbf{i})^{1/2} e^{\mathbf{i}\pi\left(\frac{ad-1}{ab}\right)u^2} \int_{-\infty}^{\infty} dt s(t) e^{\mathbf{i}\pi\frac{a}{b}\left(t-\frac{u}{\sqrt{ab}}\right)^2} \tag{20}$$

The four parameters in (20) define a characteristic matrix $\mathbf{R} = \begin{pmatrix} a & b \\ c & d \end{pmatrix}$ with $Det(\boldsymbol{R}) = 1$ performing a generalized rotation in the signal phase space. This is also followed by a special additivity property for functional compositions where $\left(O_F^{R1} \circ O_F^{R2}\right)(s(t)) = O_F^{R3}(s(t))$ where $R_3 = R_2 \cdot R_1$. This last property allows decomposing an even order $2n$ x $2n$ solution of the Yang-Baxter equation into a set of $n^2$ LCT operations. In [37], solutions for 4 x 4 matrices are given. Since the matrix elements are known one can always realize their products via decomposition into 4 blocks as $\begin{pmatrix} \mathbf{A} & \mathbf{B} \\ \mathbf{C} & \mathbf{D} \end{pmatrix}$ in which case the 4 LCTs will be independently characterized by

$$R_{11} = \mathbf{A}^2 + \mathbf{B}\cdot\mathbf{C},\ R_{12} = \mathbf{A}\cdot\mathbf{B} + \mathbf{B}\cdot\mathbf{D},\ R_{21} = \mathbf{A}\cdot\mathbf{C} + \mathbf{C}\cdot\mathbf{D},\ R_{22} = \mathbf{B}\cdot\mathbf{C} + \mathbf{D}^2 \tag{21}$$

This leads to a possible encoding of the $B_n$ algebraic structure on the phase space of some special class of signals via additional optical mixing techniques.

Notably, in [67], a complete study of all possible eigenfucntions for the LCT is given in terms of seven classes with many of them being identical with already observed "Diffractal" waveforms as first studied by Berry[68] and others [69], [70], [71]. Such waveforms are often found in the Talbot and Talbot-Lau interferometers [72], [73]. Last but not least, this instrumentation has recently been extended in material waves [74], [75], [76] leaving open the additional possibility of direct application of certain distributed computational structures in De Broglie holograms.

**5. Discussion and conclusions**

*"Truth is stranger than fiction, but it is because*
*fiction is obliged to stick to possibilities; truth isn't."*
Mark Twain, "*Following the equator: A journey around the world*"

In the previous sections, we attempted a complete arithmetization of discrete versions for certain quantum algebraic structures. The aim is towards finding methods for translating the results in a form of table lookups to be used for the transfer of the associated symmetries into signal systems and their transitions performed by appropriate filtering. Many of these structures, when analyzed with the method of global maps production offered by the combinatorial hierarchy may be resolvable as fractal sets via appropriate renormalization schemes or discrete wavelets, a theme which is still under scrutiny. Use of signal processing methods to emulate symmetries of various physical systems stands in the middle of a great dichotomy regarding the role of discrete versus continuous spaces. On the other hand, the recent discussion on finite entropy and information content of natural processes brings about the possibility of coexistence of such via special encodings. In future work, more general distributed computing systems including reservoir computing and recurrent neural network dynamics shall have to be examined in detail as possible models of hybrid machines with similar emulating capabilities.

We prefer to close this short report with a particularly intriguing analogy. Consider then the case of two engineers attempting particular signal measurements carrying different equipment. One has a spectral analyzer and the other has a heterodyne unit, possibly followed by a low pass filter. The signal to be measured has the form of a "beat" frequency or simply $\cos(\omega t) + \cos(\omega' t) = 2\cos((\omega - \omega')t/2)\cos((\omega + \omega')t/2)$. Apparently, the first engineer will be able to measure the two different frequencies $\omega, \omega'$ as the content of the signal while the second will readily isolate their half-difference to be detected as a modulated "symbol" with the half-summand as the "carrier" frequency. One cannot assert the locality of such symbol before being measured, when it remains in the form of a modulated signal envelope. This situation appears to be "contextual" enough at least in a classical frame. There is also no doubt that it represents a single, special member out of the most general category of tautologies possible.

Assume then the existence of another tautology, logical, equational or otherwise, say like $P_L(x,t) = P_R(x,t)$ each side of which encodes one of a pair of different but interlinked messages. The challenging question is whether one could find a particular member out the most general set of tautologies, possibly encoding some non-commuting algebraic structure such that an observer using a single slit optical diffraction device would read some contents of the *lhs*, while another with a double slit would read some associated contents of the *rhs*. Actually, one might conjecture such an identification of symbolically “real” objects, or signal contents as observables, to be the closest possible consistent understanding of the originally proposed “*Be-ables*” by Bell [77].

A last theme, not examined properly in this short report involves the case of severely restricted encoding and decoding maps characterized by a strong imbalance in their Kolomogorov [78] or other complexity measure thus forming what is known in cryptography as a one-way-function [79]. Combined with the advent of homomorphic encryption methods [80], [81] under which a dynamics is possible at least computationally, in a manner that allows an encrypted initial condition to remain so at all later stages of development, it could lead to the intriguing metaphor of looking at natural reality as a kind of a “*Holo-Cipher*” which by its own construction, systematically resists to any attempts of decrypting its actual, “true” contents. Some evidence towards a somewhat similar pessimistic interpretation has very recently appeared as an argument for “ugliness”, in a recent work by Hossenfelder [82]. One could perhaps detect more evidence for this in the future, as a severe exponential increase of complexity in any attempts of final unification inside a so called “Theory-of-Everything”.